\documentclass[12pt,a4paper]{article}
\usepackage{amsmath, amssymb, amsthm, amsfonts, cite, dsfont, enumerate, epsfig, float, infwarerr, mathrsfs, mathtools, mathrsfs, mathtools, relsize, stmaryrd, tabularx, txfonts}
\usepackage[normalem]{ulem}
\usepackage[colorlinks=true,citecolor=blue]{hyperref}

\usepackage{macros}
\usepackage{cite}
\usepackage{hyperref}
\usepackage[T1]{fontenc}
\usepackage{cleveref}

\hypersetup{
  colorlinks = true,
  urlcolor   = red,
  linkcolor  = blue,
  citecolor  = red
}

\makeatletter

\@addtoreset{theorem}{section}
\@addtoreset{proposition}{section}
\@addtoreset{corollary}{section}
\@addtoreset{example}{section}
\@addtoreset{remark}{section}
\@addtoreset{lemma}{section}
\@addtoreset{definition}{section}
\@addtoreset{question}{section}
\@addtoreset{conjecture}{section}

\def\@seccntformat#1{\csname the#1\endcsname\quad}
\makeatother

\makeatletter\c@MaxMatrixCols=35\makeatother   

\begin{document}
\setlength{\baselineskip}{1\baselineskip}

\setcounter{page}{1}

\title{\LARGE{An example showing that Schrijver's $\vartheta$-function need not upper bound the Shannon capacity of a graph}}

\author{Igal Sason
\thanks{
Igal Sason is with the Viterbi Faculty of Electrical and Computer Engineering and also the Department of Mathematics
(a secondary affiliation) at the Technion --- Israel Institute of Technology, Haifa 3200003, Israel. Email: eeigal@technion.ac.il.\newline
Citation: I. Sason, ``An example showing that Schrijver's $\vartheta$-function need not
   upper bound the Shannon capacity of a graph,'' {\em AIMS Mathematics}, vol.~10, no.~7,
   paper~685, pp.~15294--15301, July 2025. \url{https://www.aimspress.com/article/doi/10.3934/math.2025685}}}

\maketitle

\thispagestyle{empty}
\setcounter{page}{1}

\vspace*{-1cm}
\begin{abstract}
This letter addresses an open question concerning a variant of the Lov\'{a}sz $\vartheta$ function, which was introduced by Schrijver and
independently by McEliece et al. (1978). The question of whether this variant provides an upper bound on the Shannon capacity of
a graph was explicitly stated by Bi and Tang (2019). This letter presents an explicit example of a Tanner graph on 32 vertices,
which shows that, in contrast to the Lov\'{a}sz $\vartheta$ function, this variant does not necessarily upper bound
the Shannon capacity of a graph. The example, previously outlined by the author in a recent paper (2024), is presented here in full
detail, making it easy to follow and verify. By resolving this question, the note clarifies a subtle but significant distinction
between these two closely related graph invariants.
\end{abstract}


\section{Introduction}

The concept of the Shannon capacity of a graph, as introduced in Shannon's seminal paper (1956) on
zero-error communication \cite{Shannon56}, plays a key role in understanding the synergy and
interaction between zero-error information theory and graph theory. The zero-error capacity
of a discrete memoryless channel (DMC) equals the Shannon capacity of the corresponding confusability graph,
whose vertices represent the input symbols of the channel, and any pair of its vertices are adjacent
if they represent input symbols that can be confused by the channel (i.e., conditioned on each of these
two input symbols, an identical output symbol can be received with some positive probabilities).
The significance of the Shannon capacity of graphs, and the hardness of its computability in general,
are highlighted in various survey papers \cite{Alon02,Alon19,Jurkiewicz14,KornerO98}.

The aim of this note is to present an example demonstrating that, in contrast to the Lov\'{a}sz
$\vartheta$-function, which upper bounds the Shannon capacity of a graph \cite{Lovasz79_IT}, the variant
introduced by Schrijver does not \cite[Example~5.24]{Sason24}. This resolves a query concerning the variant
proposed by Schrijver, which is identical to the function independently presented by McEliece {\em et al.} (1978)
\cite{McElieceRR78,Schrijver79}, and that was posed as an open question by Bi and Tang \cite{BiT19}. They wrote in \cite{BiT19}:
“In fact, it is quite difficult to disprove that the Schrijver $\vartheta$-function is an upper bound, because
at least for graphs $\Gr{G}$ of moderate size the two values $\vartheta(\Gr{G})$ (the Lov\'{a}sz $\vartheta$-function
of $\Gr{G}$) and $\vartheta'(\Gr{G})$ (the Schrijver $\vartheta$-function of $\Gr{G}$) are very close to each other.”
Such an example is presented in Section~\ref{section: example}, following a brief overview of preliminary material in
Section~\ref{section: preliminaries}.

\section{Preliminaries}
\label{section: preliminaries}

In the following, $\indnum{\Gr{G}}$, $\vartheta(\Gr{G})$, $\vartheta'(\Gr{G})$, and $\Theta(\Gr{G})$ respectively denote the independence
number, Lov\'{a}sz $\vartheta$-function, Schrijver's $\vartheta$-function, and the Shannon capacity of a simple, finite, and undirected
graph~$\Gr{G}$.
We refer the reader to our recent paper (see \cite[Section~2]{Sason24}) for preliminaries, definitions, and an account of the properties of the
Lov\'{a}sz and Schrijver's $\vartheta$-functions of a graph, as well as the Shannon capacity of a graph.
In order to introduce the Shannon capacity of a graph, we need the notion of a {\em strong product} of graphs.
\begin{definition}[Strong products and strong powers of graphs]
\label{def:strong product of graphs}
{\em Let $\Gr{G}$ and $\Gr{H}$ be graphs.
The strong product $\Gr{G} \boxtimes \Gr{H}$ is a graph with
a vertex set $\V{\Gr{G} \boxtimes \Gr{H}} = \V{\Gr{G}} \times \V{\Gr{H}}$ (the Cartesian product),
and distinct vertices $(g, h)$ and $(g', h')$ are adjacent in $\Gr{G} \boxtimes \Gr{H}$ if
and only if one of the following three conditions holds:
(1) $g = g'$ and $\{h, h'\} \in \E{\Gr{H}}$, (2) $\{g, g'\} \in \E{\Gr{G}}$ and $h = h'$, or
(3) $\{g, g'\} \in \E{\Gr{G}}$ and $\{h, h'\} \in \E{\Gr{H}}$.
\newline
Strong products are therefore commutative and associative up to graph isomorphisms. The {\em $k$-fold
strong power} of $\Gr{G}$ is given by $\Gr{G}^{\boxtimes \, k} \eqdef \Gr{G} \boxtimes \cdots \boxtimes \Gr{G}$,
where $\Gr{G}$ is multiplied by itself $k$ times.}
\end{definition}

The following three results are used in this note.
\begin{theorem} \cite{Shannon56,Lovasz79_IT}
\label{theorem: Shannon capacity}
{\em For every simple graph $\Gr{G}$,
\begin{align}
\Theta(\Gr{G}) &\triangleq \sup_{k \in \naturals} \sqrt[k]{\indnum{\Gr{G}^{\boxtimes \, k}}} \label{eq1: capacity} \\
&= \lim_{k \to \infty} \sqrt[k]{\indnum{\Gr{G}^{\boxtimes \, k}}},  \label{eq2: capacity}
\end{align}
and
\begin{align}
\indnum{\Gr{G}} \leq \Theta(\Gr{G}) \leq \vartheta(\Gr{G}).
\end{align}}
\end{theorem}

\begin{theorem} \cite{Lovasz79_IT,Schrijver79}
\label{theorem: theta functions}
{\em For every graph $\Gr{G}$,
\begin{align}
\label{eq2:17.11.23}
\indnum{\Gr{G}} \leq \vartheta'(\Gr{G}) \leq \vartheta(\Gr{G}).
\end{align}}
\end{theorem}

\begin{theorem} \cite{Lovasz79_IT}
\label{theorem: vartheta et}
{\em Let $\Gr{G}$ be $d$-regular of order $n$, and let $\Eigval{n}{\Gr{G}}$ be the smallest eigenvalue of its adjacency matrix. Then,
\begin{align}
\label{eq: Lovasz79 - Theorem 9}
\vartheta(\Gr{G}) \leq -\frac{n \, \Eigval{n}{\Gr{G}}}{d - \Eigval{n}{\Gr{G}}},
\end{align}
with equality if $\Gr{G}$ is an edge-transitive graph.}
\end{theorem}

The interested reader is also referred to a recent survey paper \cite{Zhou25} on upper bounds on the independence number of a graph,
including the Lov\'{a}sz and Schrijver's $\vartheta$ functions (see \eqref{eq2:17.11.23}).

In light of Theorems~\ref{theorem: Shannon capacity} and~\ref{theorem: theta functions}, the following question was left open since \cite{Schrijver79},
and recently resolved in \cite[Example~5.24]{Sason24}.

\begin{question}
\label{question: Shannon capacity}
{\em Could the bound on the Shannon capacity,
$\Theta(\Gr{G}) \leq \vartheta(\Gr{G})$,
be improved by the bound
$\Theta(\Gr{G}) \leq \vartheta'(\Gr{G})$?}
\end{question}
We give a negative answer to Question~\ref{question: Shannon capacity} by
providing a complete and detailed presentation of the example in \cite[Example~5.24]{Sason24}. The main result in this note
is thus an explicit and detailed example demonstrating that
\begin{align}
\label{eq:main result}
\Theta(\Gr{G}) \not\leq \vartheta'(\Gr{G}).
\end{align}

\section{Example}
\label{section: example}

Let $\Gr{G}$ be the Gilbert graph on~32 vertices, where
\begin{align}
\label{eq: Gilbert graph}
\V{\Gr{G}} = \{0,1\}^5, \qquad \E{\Gr{G}} = \Bigl\{ \underline{u}, \underline{v} \in \{0,1\}^5: \; 1 \leq d_{\mathrm{H}}(\underline{u}, \underline{v}) \leq 2 \Bigr\},
\end{align}
so, every two vertices are adjacent if and only if the Hamming distance between their corresponding binary 5-length tuples is either~1 or 2. Label
each vertex in $\{0,1\}^5$ by its corresponding decimal value in $\{0, 1, \ldots, 31\}$, and let the $i$-th row and column of the adjacency matrix ${\bf{A}}$
correspond to the vertex labeled $i-1$, for each $i \in \OneTo{32}$. Then, the following holds:
\begin{itemize}
\item The graph $\Gr{G}$ is 15-regular, since $\binom{5}{1}+\binom{5}{2} = 5 + 10 = 15$.
\item To motivate why this graph is a suitable candidate for such an example, it is instructive to observe the following structural properties:
\begin{enumerate}
\item $\Gr{G}$ is vertex-transitive, edge-transitive, and distance-regular, but not strongly regular.
\item The complement graph $\CGr{G}$ is 16-regular and vertex-transitive, but it is not edge-transitive nor distance-regular (and thus not strongly regular).
\end{enumerate}
These properties, which can be observed using the SageMath software~\cite{SageMath}, are not essential for verifying
the subsequent steps of the example. However, if the complement graph $\CGr{G}$ were edge-transitive in addition to being vertex-transitive,
or if it were strongly regular, then by combining Theorem~5.9 and Eq.~(2.43) of \cite{Sason24} with Proposition~1 of \cite{Sason23},
it would follow that $\vartheta(\Gr{G}) = \vartheta'(\Gr{G})$, and hence $\Gr{G}$ could not serve as a counterexample to Question~\ref{question: Shannon capacity}
(since $\Theta(\Gr{G}) \leq \vartheta(\Gr{G})$ by Theorem~1 in \cite{Lovasz79_IT}).
\item The independence number of $\Gr{G}$ is $\indnum{\Gr{G}} = 4$. An example of such a maximal independent set of $\Gr{G}$:
\begin{align}
\label{maximal independent set of size 4}
\bigl\{(1, 0, 0, 1, 0), \; (0, 1, 1, 1, 0), \; (0, 0, 0, 0, 1), \; (1, 1, 1, 0, 1)\bigr\}.
\end{align}
\item Solving the following minimization problem for $\vartheta'(\Gr{G})$ \cite{Schrijver79} (see also \cite{Zhou25}):
\vspace*{0.1cm}
\begin{eqnarray}
\label{eq: SDP problem - Schrijver theta-function}
\mbox{\fbox{$
\begin{array}{l}
\text{minimize} \; \; \lambda_{\max}({\bf{X}})  \\
\text{subject to} \\
\begin{cases}
{\bf{X}} \in \mathcal{S}_{32} \\
A_{i,j} = 0  \; \Rightarrow \;  X_{i,j} \geq 1, \quad \forall \, i,j \in \{1, \ldots, 32\}
\end{cases}
\end{array}$}}
\end{eqnarray}
where $\mathcal{S}_{32}$ denotes the set of all the $32 \times 32$ real symmetric matrices,
is a dual semidefinite programming (SDP) problem with strong duality.
We aim to show that
$\vartheta'(\Gr{G}) = 4 = \indnum{\Gr{G}}$, relying on the CVX software \cite{CVX}.
\item To that end, consider the feasible solution of \eqref{eq: SDP problem - Schrijver theta-function}
$${\bf{X}} = [{\bf{X}}_{\mathrm{l}}, \; {\bf{X}}_{\mathrm{r}}],$$
where ${\bf{X}}_{\mathrm{l}}$ and ${\bf{X}}_{\mathrm{r}}$ are given by the following $32 \times 16$ submatrices:
\begin{align}
\label{eq:X_l}
{\bf{X}}_{\mathrm{l}} =
\begin{pmatrix}
+1 & -1 &  -1 &  -1 &  -1 &  -1 &  -1 &  +1 &  -1 &  -1 &  -1 &  +1 &  -1 &  +1 &  +1 &  +1 \\
-1 & +1 &  -1 &  -1 &  -1 &  -1 &  +1 &  -1 &  -1 &  -1 &  +1 &  -1 &  +1 &  -1 &  +1 &  +1 \\
-1 & -1 &  +1 &  -1 &  -1 &  +1 &  -1 &  -1 &  -1 &  +1 &  -1 &  -1 &  +1 &  +1 &  -1 &  +1 \\
-1 & -1 &  -1 &  +1 &  +1 &  -1 &  -1 &  -1 &  +1 &  -1 &  -1 &  -1 &  +1 &  +1 &  +1 &  -1 \\
-1 & -1 &  -1 &  +1 &  +1 &  -1 &  -1 &  -1 &  -1 &  +1 &  +1 &  +1 &  -1 &  -1 &  -1 &  +1 \\
-1 & -1 &  +1 &  -1 &  -1 &  +1 &  -1 &  -1 &  +1 &  -1 &  +1 &  +1 &  -1 &  -1 &  +1 &  -1 \\
-1 & +1 &  -1 &  -1 &  -1 &  -1 &  +1 &  -1 &  +1 &  +1 &  -1 &  +1 &  -1 &  +1 &  -1 &  -1 \\
+1 & -1 &  -1 &  -1 &  -1 &  -1 &  -1 &  +1 &  +1 &  +1 &  +1 &  -1 &  +1 &  -1 &  -1 &  -1 \\
-1 & -1 &  -1 &  +1 &  -1 &  +1 &  +1 &  +1 &  +1 &  -1 &  -1 &  -1 &  -1 &  -1 &  -1 &  +1 \\
-1 & -1 &  +1 &  -1 &  +1 &  -1 &  +1 &  +1 &  -1 &  +1 &  -1 &  -1 &  -1 &  -1 &  +1 &  -1 \\
-1 & +1 &  -1 &  -1 &  +1 &  +1 &  -1 &  +1 &  -1 &  -1 &  +1 &  -1 &  -1 &  +1 &  -1 &  -1 \\
+1 & -1 &  -1 &  -1 &  +1 &  +1 &  +1 &  -1 &  -1 &  -1 &  -1 &  +1 &  +1 &  -1 &  -1 &  -1 \\
-1 & +1 &  +1 &  +1 &  -1 &  -1 &  -1 &  +1 &  -1 &  -1 &  -1 &  +1 &  +1 &  -1 &  -1 &  -1 \\
+1 & -1 &  +1 &  +1 &  -1 &  -1 &  +1 &  -1 &  -1 &  -1 &  +1 &  -1 &  -1 &  +1 &  -1 &  -1 \\
+1 & +1 &  -1 &  +1 &  -1 &  +1 &  -1 &  -1 &  -1 &  +1 &  -1 &  -1 &  -1 &  -1 &  +1 &  -1 \\
+1 & +1 &  +1 &  -1 &  +1 &  -1 &  -1 &  -1 &  +1 &  -1 &  -1 &  -1 &  -1 &  -1 &  -1 &  +1 \\
-1 & -1 &  -1 &  +1 &  -1 &  +1 &  +1 &  +1 &  -1 &  +1 &  +1 &  +1 &  +1 &  +1 &  +1 &  {\bf{+3}} \\
-1 & -1 &  +1 &  -1 &  +1 &  -1 &  +1 &  +1 &  +1 &  -1 &  +1 &  +1 &  +1 &  +1 &  {\bf{+3}} &  +1 \\
-1 & +1 &  -1 &  -1 &  +1 &  +1 &  -1 &  +1 &  +1 &  +1 &  -1 &  +1 &  +1 &  {\bf{+3}} &  +1 &  +1 \\
+1 & -1 &  -1 &  -1 &  +1 &  +1 &  +1 &  -1 &  +1 &  +1 &  +1 &  -1 &  {\bf{+3}} &  +1 &  +1 &  +1 \\
-1 & +1 &  +1 &  +1 &  -1 &  -1 &  -1 &  +1 &  +1 &  +1 &  +1 &  {\bf{+3}} &  -1 &  +1 &  +1 &  +1 \\
+1 & -1 &  +1 &  +1 &  -1 &  -1 &  +1 &  -1 &  +1 &  +1 &  {\bf{+3}} &  +1 &  +1 &  -1 &  +1 &  +1 \\
+1 & +1 &  -1 &  +1 &  -1 &  +1 &  -1 &  -1 &  +1 &  {\bf{+3}} &  +1 &  +1 &  +1 &  +1 &  -1 &  +1 \\
+1 & +1 &  +1 &  -1 &  +1 &  -1 &  -1 &  -1 &  {\bf{+3}} &  +1 &  +1 &  +1 &  +1 &  +1 &  +1 &  -1 \\
-1 & +1 &  +1 &  +1 &  +1 &  +1 &  +1 &  {\bf{+3}} &  -1 &  -1 &  -1 &  +1 &  -1 &  +1 &  +1 &  +1 \\
+1 & -1 &  +1 &  +1 &  +1 &  +1 &  {\bf{+3}} &  +1 &  -1 &  -1 &  +1 &  -1 &  +1 &  -1 &  +1 &  +1 \\
+1 & +1 &  -1 &  +1 &  +1 &  {\bf{+3}} &  +1 &  +1 &  -1 &  +1 &  -1 &  -1 &  +1 &  +1 &  -1 &  +1 \\
+1 & +1 &  +1 &  -1 &  {\bf{+3}} &  +1 &  +1 &  +1 &  +1 &  -1 &  -1 &  -1 &  +1 &  +1 &  +1 &  -1 \\
+1 & +1 &  +1 &  {\bf{+3}} &  -1 &  +1 &  +1 &  +1 &  -1 &  +1 &  +1 &  +1 &  -1 &  -1 &  -1 &  +1 \\
+1 & +1 &  {\bf{+3}} &  +1 &  +1 &  -1 &  +1 &  +1 &  -1 &  -1 &  +1 &  +1 &  -1 &  -1 &  +1 &  -1 \\
+1 & {\bf{+3}} &  +1 &  +1 &  +1 &  +1 &  -1 &  +1 &  +1 &  +1 &  -1 &  +1 &  -1 &  +1 &  -1 & -1 \\
{\bf{+3}} & +1 &  +1 &  +1 &  +1 &  +1 &  +1 &  -1 &  +1 &  +1 &  +1 &  -1 &  +1 &  -1 &  -1 & -1
\end{pmatrix},
\end{align}
and
\begin{align}
\label{eq:X_r}
{\bf{X}}_{\mathrm{r}} =
\begin{pmatrix}
-1  &  -1  &  -1  &  +1  &  -1  &   +1  &   +1  &   +1  &  -1  &  +1  &   +1  &   +1  &   +1  &   +1  &   +1  &   {\bf{+3}}  \\
-1  &  -1  &  +1  &  -1  &  +1  &   -1  &   +1  &   +1  &  +1  &  -1  &   +1  &   +1  &   +1  &   +1  &   {\bf{+3}}  &   +1  \\
-1  &  +1  &  -1  &  -1  &  +1  &   +1  &  -1   &   +1  &  +1  &  +1  &   -1  &   +1  &   +1  &   {\bf{+3}}  &   +1  &   +1  \\
+1  &  -1  &  -1  &  -1  &  +1  &   +1  &  +1   &   -1  &  +1  &  +1  &   +1  &   -1  &   {\bf{+3}}  &   +1  &   +1  &   +1  \\
-1  &  +1  &  +1  &  +1  &  -1  &   -1  &  -1   &   +1  &  +1  &  +1  &   +1  &   {\bf{+3}}  &   -1  &   +1  &   +1  &   +1  \\
+1  &  -1  &  +1  &  +1  &  -1  &   -1  &  +1   &   -1  &  +1  &  +1  &   {\bf{+3}}  &   +1  &   +1  &   -1  &   +1  &   +1  \\
+1  &  +1  &  -1  &  +1  &  -1  &   +1  &  -1   &   -1  &  +1  &  {\bf{+3}}  &   +1  &   +1  &   +1  &   +1  &   -1  &   +1  \\
+1  &  +1  &  +1  &  -1  &  +1  &   -1  &  -1   &   -1  &  {\bf{+3}}  &  +1  &   +1  &   +1  &   +1  &   +1  &   +1  &   -1  \\
-1  &  +1  &  +1  &  +1  &  +1  &   +1  &  +1   &   {\bf{+3}}  &  -1  &  -1  &  -1   &   +1  &   -1  &   +1  &   +1  &   +1  \\
+1  &  -1  &  +1  &  +1  &  +1  &   +1  &  {\bf{+3}}   &   +1  &  -1  &  -1  &   1   &   -1  &   +1  &   -1  &   +1  &   +1  \\
+1  &  +1  &  -1  &  +1  &  +1  &   {\bf{+3}}  &  +1   &   +1  &  -1  &  +1  &  -1   &   -1  &   +1  &   +1  &   -1  &   1   \\
+1  &  +1  &  +1  &  -1  &  {\bf{+3}}  &   +1  &  +1   &   +1  &  +1  &  -1  &  -1   &   -1  &   +1  &   +1  &   +1  &  -1   \\
+1  &  +1  &  +1  &  {\bf{+3}}  &  -1  &   +1  &  +1   &   +1  &  -1  &  +1  &  +1   &   +1  &   -1  &   -1  &   -1  &   1   \\
+1  &  +1  &  {\bf{+3}}  &  +1  &  +1  &   -1  &  +1   &   +1  &  +1  &  -1  &  +1   &   +1  &   -1  &   -1  &   +1  &  -1   \\
+1  &  {\bf{+3}}  &  +1  &  +1  &  +1  &   +1  &  -1   &   +1  &  +1  &  +1  &  -1   &   +1  &   -1  &   +1  &   -1  &  -1   \\
{\bf{+3}}  &  +1  &  +1  &  +1  &  +1  &   +1  &  +1   &   -1  &  +1  &  +1  &  +1   &   -1  &   +1  &   -1  &   -1  &  -1   \\
+1  &  -1  &  -1  &  -1  &  -1  &   -1  &  -1   &   +1  &  -1  &  -1  &  -1   &   +1  &   -1  &   +1  &   +1  &  +1   \\
-1  &  +1  &  -1  &  -1  &  -1  &   -1  &  +1   &   -1  &  -1  &  -1  &  +1   &   -1  &   +1  &   -1  &   +1  &  +1   \\
-1  &  -1  &  +1  &  -1  &  -1  &   +1  &  -1   &   -1  &  -1  &  +1  &  -1   &   -1  &   +1  &   +1  &   -1  &  +1   \\
-1  &  -1  &  -1  &  +1  &  +1  &   -1  &  -1   &   -1  &  +1  &  -1  &  -1   &   -1  &   +1  &   +1  &   +1  &  -1   \\
-1  &  -1  &  -1  &  +1  &  +1  &   -1  &  -1   &   -1  &  -1  &  +1  &  +1   &   +1  &   -1  &   -1  &   -1  &  +1   \\
-1  &  -1  &  +1  &  -1  &  -1  &   +1  &  -1   &   -1  &  +1  &  -1  &  +1   &   +1  &   -1  &   -1  &   +1  &  -1   \\
-1  &  +1  &  -1  &  -1  &  -1  &   -1  &  +1   &   -1  &  +1  &  +1  &  -1   &   +1  &   -1  &   +1  &   -1  &  -1   \\
+1  &  -1  &  -1  &  -1  &  -1  &   -1  &  -1   &   +1  &  +1  &  +1  &  +1   &   -1  &   +1  &   -1  &   -1  &  -1   \\
-1  &  -1  &  -1  &  +1  &  -1  &   +1  &  +1   &   +1  &  +1  &  -1  &  -1   &   -1  &   -1  &   -1  &   -1  &  +1   \\
-1  &  -1  &  +1  &  -1  &  +1  &   -1  &  +1   &   +1  &  -1  &  +1  &  -1   &   -1  &   -1  &   -1  &   +1  &  -1   \\
-1  &  +1  &  -1  &  -1  &  +1  &   +1  &  -1   &   +1  &  -1  &  -1  &  +1   &   -1  &   -1  &   +1  &   -1  &  -1   \\
+1  &  -1  &  -1  &  -1  &  +1  &   +1  &  +1   &   -1  &  -1  &  -1  &  -1   &   +1  &   +1  &   -1  &   -1  &  -1   \\
-1  &  +1  &  +1  &  +1  &  -1  &   -1  &  -1   &   +1  &  -1  &  -1  &  -1   &   +1  &   +1  &   -1  &   -1  &  -1   \\
+1  &  -1  &  +1  &  +1  &  -1  &   -1  &  +1   &   -1  &  -1  &  -1  &  +1   &   -1  &   -1  &   +1  &   -1  &  -1   \\
+1  &  +1  &  -1  &  +1  &  -1  &   +1  &  -1   &   -1  &  -1  &  +1  &  -1   &   -1  &   -1  &   -1  &   +1  &  -1   \\
+1  &  +1  &  +1  &  -1  &  +1  &   -1  &  -1   &   -1  &  +1  &  -1  &  -1   &   -1  &   -1  &   -1  &   -1  &   1
\end{pmatrix},
\end{align}
which implies that the largest eigenvalue of the $32 \times 32$ matrix ${\bf{X}} = [{\bf{X}}_{\mathrm{l}}, \; {\bf{X}}_{\mathrm{r}}]$
is equal to
\begin{align}
\label{eq: max eig}
\lambda_{\max}({\bf{X}}) = 4.
\end{align}
\item Combining the left side of \eqref{eq2:17.11.23} with \eqref{maximal independent set of size 4},
\eqref{eq: SDP problem - Schrijver theta-function}, and \eqref{eq: max eig}
yields $4 = \indnum{\Gr{G}} \leq \vartheta'(\Gr{G}) \leq \lambda_{\max}({\bf{X}}) = 4,$ which gives
\begin{align}
\label{eq: schrijver's theta and independence number}
\vartheta'(\Gr{G}) = 4 = \indnum{\Gr{G}}.
\end{align}
Note that all the entries of the matrix ${\bf{X}} = [{\bf{X}}_{\mathrm{l}}, \; {\bf{X}}_{\mathrm{r}}]$, as given in \eqref{eq:X_l}
and \eqref{eq:X_r}, take the values $3, 1, -1$. Furthermore, $X_{i,j}=3$ if and only if $i+j=33$ and $i,j \in \OneTo{32}$ (i.e.,
the bolded entries in \eqref{eq:X_l} and \eqref{eq:X_r} take the value~3, and they form the antidiagonal of ${\bf{X}}$).
\item Although not needed directly for establishing the counterexample, the graph $\Gr{G}$ is 15-regular and edge-transitive on 32 vertices,
with $\lambda_{\min}(\Gr{G}) = -3$. Hence, by Theorem~\ref{theorem: vartheta et}, the Lov\'{a}sz $\vartheta$-function of $\Gr{G}$ is given by
\begin{align}
\label{eq: Lovasz number}
\vartheta(\Gr{G}) = -\frac{n \lambda_{\min}(\Gr{G})}{d(\Gr{G}) - \lambda_{\min}(\Gr{G})} = \frac{32 \cdot 3}{15+3} = 5 \frac13.
\end{align}
This can also be verified numerically by solving the dual SDP problem
\vspace*{0.1cm}
\begin{eqnarray}
\label{eq: SDP problem - Lovasz theta-function}
\mbox{\fbox{$
\begin{array}{l}
\text{minimize} \; \; \lambda_{\max}({\bf{X}})  \\
\text{subject to} \\
\begin{cases}
{\bf{X}} \in \mathcal{S}_{32} \\
A_{i,j} = 0  \; \Rightarrow \;  X_{i,j} = 1, \quad \forall \, i,j \in \{1, \ldots, 32\},
\end{cases}
\end{array}$}}
\end{eqnarray}
differing from the dual SDP problem in \eqref{eq: SDP problem - Schrijver theta-function} in the stronger requirement
that, whenever $A_{i,j} = 0$, the entry $X_{i,j}$ is restricted to be equal to~1, rather than imposing the weaker condition that $X_{i,j} \geq 1$.
It is interesting to note the following relation: let ${\bf{X}}$ and ${\bf{\hat{X}}}$ be the matrices that solve the SDP problems in
\eqref{eq: SDP problem - Schrijver theta-function} and \eqref{eq: SDP problem - Lovasz theta-function}, respectively. Then, the following holds:
\begin{enumerate}[(1)]
\item If $X_{i,j} = -1$, then $\hat{X}_{i,j} = -\tfrac79$;
\item If $X_{i,j} \in \{1, 3\}$, then $\hat{X}_{i,j} = 1$;
\item All entries in the antidiagonal of ${\bf{X}}$ are equal to~3, and no other entry in ${\bf{X}}$ is equal to~3.
\end{enumerate}
This transforms the optimized matrix ${\bf{X}}$ in \eqref{eq: SDP problem - Schrijver theta-function} to the one in \eqref{eq: SDP problem - Lovasz theta-function}
and vice versa.
Note that combining with \eqref{eq: schrijver's theta and independence number}, we have for the graph $\Gr{G}$,
\begin{align}
\label{eq: ranking}
4 = \indnum{\Gr{G}} = \vartheta'(\Gr{G}) < \vartheta(\Gr{G}) = 5 \frac13,
\end{align}
so $\vartheta'(\Gr{G})$ coincides with the independence number of $\Gr{G}$, and it is strictly smaller than $\vartheta(\Gr{G})$.
\item
By the SageMath software, we have
$\indnum{\Gr{G} \boxtimes \Gr{G}} = 20$,
and the strong product graph $\Gr{G} \boxtimes \Gr{G}$ has $368,640$ such maximal independent sets of size 20.
In fact, it is enough for our purpose to establish the weaker result where
\begin{align}
\label{eq: independence number of strong product}
\indnum{\Gr{G} \boxtimes \Gr{G}} \geq 20.
\end{align}
To that end, an example of a (maximal) independent set of size~20 for the strong product of $\Gr{G}$ by itself,
$\Gr{G} \boxtimes \Gr{G}$, is given by
\begin{align}
& \Bigl\{ ((1, 1, 0, 0, 0), (1, 1, 1, 1, 1)), \quad  ((1, 0, 1, 0, 0), (1, 1, 0, 0, 0)), \quad ((0, 1, 1, 0, 0), (0, 0, 1, 1, 0)), \nonumber \\
& ((1, 1, 1, 0, 0), (0, 0, 0, 0, 1)), \quad ((1, 0, 0, 1, 0), (0, 0, 1, 0, 1)), \quad  ((0, 1, 0, 1, 0), (1, 0, 0, 0, 0)), \nonumber \\
& ((1, 1, 0, 1, 0), (0, 1, 0, 1, 0)), \quad  ((0, 0, 1, 1, 0), (0, 1, 0, 1, 1)), \quad  ((1, 0, 1, 1, 0), (1, 0, 1, 1, 0)), \nonumber \\
& ((0, 1, 1, 1, 0), (1, 1, 1, 0, 1)), \quad ((1, 0, 0, 0, 1), (0, 0, 0, 1, 0)),  \quad  ((0, 1, 0, 0, 1), (0, 1, 0, 0, 1)), \nonumber \\
&  ((1, 1, 0, 0, 1), (1, 0, 1, 0, 0)), \quad  ((0, 0, 1, 0, 1), (1, 0, 1, 0, 1)), \quad  ((1, 0, 1, 0, 1), (0, 1, 1, 1, 1)), \nonumber \\
&  ((0, 1, 1, 0, 1), (1, 1, 0, 1, 0)), \quad  ((0, 0, 0, 1, 1), (1, 1, 1, 1, 0)), \quad  ((1, 0, 0, 1, 1), (1, 1, 0, 0, 1)), \nonumber \\
&  ((0, 1, 0, 1, 1), (0, 0, 1, 1, 1)), \quad ((0, 0, 1, 1, 1), (0, 0, 0, 0, 0)) \Bigr\}.
\end{align}
\item Consequently, it follows from \eqref{eq1: capacity}, \eqref{eq: schrijver's theta and independence number}, and
\eqref{eq: independence number of strong product} that
\begin{align}
\Theta(\Gr{G}) \geq \sqrt{\indnum{\Gr{G} \boxtimes \Gr{G}}} \geq \sqrt{20} > 4 = \vartheta'(\Gr{G}),
\end{align}
demonstrating that $\vartheta'(\Gr{G})$ does not serve as a universal upper bound on the Shannon capacity $\Theta(\Gr{G})$,
in contrast to the Lov\'{a}sz $\vartheta$-function. On the other hand, $\vartheta'(\Gr{G})$ yields an improved upper bound on the
independence number $\indnum{\Gr{G}}$ than $\vartheta(\Gr{G})$, and in fact, $\vartheta'(\Gr{G})$ provides a tight upper bound
for the considered graph $\Gr{G}$ in \eqref{eq: Gilbert graph} (see \eqref{eq: ranking}).
\end{itemize}

\end{document}